\documentclass[a4paper,12pt]{amsart}
\usepackage{amsmath,amsthm,amssymb, amscd, amsfonts, color, pst-node}
\usepackage[arrow,matrix]{xy}

\theoremstyle{plain}
\newtheorem{theorem}{Theorem}[section]
\newtheorem{proposition}[theorem]{Proposition}

\newtheorem{corollary}[theorem]{Corollary}

\newtheorem{problem}[theorem]{Problem}

\theoremstyle{plain} \numberwithin{equation}{section}
\theoremstyle{definition}

\newtheorem{definition}[theorem]{Definition}
\newtheorem{remark}[theorem]{Remark}

\newtheorem{example}[theorem]{Example}

\newtheorem*{CRPT}{Cohomological rigidity problem for toric manifolds}

\newtheorem*{HRPRT}{Homotopical rigidity problem for real toric manifolds}
\newtheorem*{CRPTT}{Cohomological rigidity problem for topological toric manifolds}

\newtheorem*{SCRPT}{Strong cohomological rigidity problem for toric manifolds}

\newcommand{\Z}{\mathbb{Z}}
\newcommand{\N}{\mathbb{N}}
\newcommand{\Q}{\mathbb{Q}}

\newcommand{\C}{\mathbb{C}}
\newcommand{\CP}{\mathbb{C}P}

\newcommand{\R}{\mathbb{R}}

\newcommand{\K}{\mathcal{K}}

\def\vput(#1)#2{\cnode*(#1){1pt}{#2}}

\def\k{\mathbf{k}}

\DeclareMathOperator{\Ker}{Ker}
\DeclareMathOperator{\Hom}{Hom}

\def\Diff#1{\mathcal D_#1}
\def\Ori#1{\mathcal O_#1}
\def\Sym#1{\mathcal S_#1}

\begin{document}
\title[Rigidity problems in toric topology]{Rigidity problems in toric topology, a survey}

\author[S. Choi]{Suyoung Choi}
\address{Department of Mathematics, Ajou University, San 5, Woncheon-dong, Yeongtong-gu, Suwon 443-749, Korea}
\email{schoi@ajou.ac.kr}

\author[M. Masuda]{Mikiya Masuda}
\address{Department of Mathematics, Osaka City University, Sugimoto, Sumiyoshi-ku, Osaka 558-8585, Japan}
\email{masuda@sci.osaka-cu.ac.jp}
\author[D. Y. Suh]{Dong Youp Suh}
\address{Department of Mathematical Sciences, KAIST, 335 Gwahangno, Yu-sung Gu, Daejeon 305-701, Korea}
\email{dysuh@math.kaist.ac.kr}
\thanks{The first author was supported by the Japanese Society for the Promotion of Sciences (JSPS grant no. P09023)
and the second author was partially supported by Grant-in-Aid for Scientific Research 22540094. The third author was partially supported by Basic Science Research Program through the National Research Foundation of Korea(NRF) funded by the Ministry of Education, Science and Technology (2010-0001651)}


\begin{abstract}
Several rigidity problems in toric topology are addressed in \cite{ma-su08}.  In this paper, we survey results on those problems including recent development.
\end{abstract}

\date{\today}
\maketitle

\section{Introduction}
As is well-known, cohomology ring does not distinguish closed smooth manifolds up to diffeomorphism or homeomorphism in general.  However, it does if we restrict our attention to a reasonably small class of objects.
For instance, simply connected closed smooth $4$-manifolds are classified up to homeomorphism by Freedman \cite{freedman:1982} using their integral cohomology rings.

In the proceedings of 2006 Osaka conference on toric topology, the second and third authors \cite{ma-su08} proposed a rather bold question asking whether the class of toric manifolds is topologically determined by their integral cohomology rings. In fact, the cohomological rigidity problem for toric manifolds asks whether two toric manifolds are homeomorphic or even diffeomorphic if their integral cohomology rings are isomorphic as graded rings.

The affirmative solution to the problem seems implausible at first glance, but no example providing the negative solution to the problem has been found yet. Instead, during the last five years, many results supporting the affirmative solution to the problem have appeared.  In order to solve the problem, we either look for an example for the negative solution, or prove for the positive solution to the problem. To the end of the latter direction, the class of toric manifolds seems too large to handle at once. Therefore, it is reasonable to focus further down to a smaller subclass of toric manifolds to get any affirmative results.

So far, all affirmative results are on (generalized) Bott manifolds, and these results are surveyed in Section~\ref{section1}. See Section~\ref{section1} for the definition of (generalized) Bott manifolds. Even though there are several results on cohomological rigidity for (generalized) Bott manifolds, it is still open, and it will be interesting to settle the problem for them.

The notion of quasitoric manifold is introduced in \cite{da-ja91} as a topological analogue of toric manifolds. See Section~\ref{section2}
for the definition of quasitoric manifold. In particular, all projective toric manifolds are quasitoric manifolds, but the converse is not always true. A similar cohomological rigidity problem can be raised for quasitoric manifolds, and so far, all the affirmative results are found for quasitoric manifolds over products of simplices, i.e., quasitoric manifolds whose orbit spaces are products of simplices. Note that a (generalized) Bott manifold is an example of a quasitoric manifold over a product of simplices.

Quasitoric manifold is not the only topological analogue of toric manifold. In fact, topological toric manifold is defined in \cite{is-fu-ma-pre}, and, in particular, all toric manifolds and quasitoric manifolds are topological toric manifolds. Therefore, topological toric manifold can be considered as the correct topological analogue of toric manifold. A cohomological rigidity problem for topological toric manifolds can be raised similarly, but so far there is no results on this problem except for those for toric or quasitoric manifolds discussed above. Torus manifold defined in \cite{ha-ma03} and \cite{masu99} is a further generalization of toric manifold, however, the class of torus manifolds is too large to be cohomologically rigid. The developments on rigidity of quasitoric, topological toric, and torus manifolds are surveyed in Section~\ref{section2}. We also mention some rigidity result for symplectic toric manifolds there.

In Section~\ref{section3}, real versions of toric manifold,  quasitoric manifold called small cover, and their corresponding cohomological rigidity problems are discussed. In these cases, however, integral cohomologies are quite difficult to calculate, so the rigidity problems with integral cohomologies are unrealistic. Instead, we consider cohomology with $\Z/2$-coefficients to get any reasonable results. Notable results in this case are the complete settlement of cohomological rigidity for real (generalized) Bott manifolds in \cite{ka-ma09} and \cite{masu10}.
Namely, the problem is affirmative for real Bott manifolds, but negative for generalized real Bott manifolds.

The strong cohomological rigidity problem for toric manifolds asks whether an isomorphism between the cohomology rings of two toric manifolds
is induced from a homeomorphism or a diffeomorphism between the manifolds.  The strong cohomological rigidity problem for quasitoric or topological toric manifolds can be asked similarly, but note that this problem with diffeomorphism has the negative answer by the work of Friedman and Morgan \cite{Fr-Mo-1988}. There are some positive results for the problem with diffeomorphism, and these together with some related results are surveyed in Section~\ref{section4}.

Toric manifolds and quasitoric manifolds are associated with combinatorial objects which are simplicial complexes and simple convex polytopes, respectively. Even though there are some counterexamples, often such combinatorial objects are completely determined by the cohomology rings of the associated manifolds, and in this case such simplicial complexes or polytopes are said to be (cohomologically) rigid. For such combinatorial objects, the rigidity problem  asks to determine (cohomologically) rigid complexes or polytopes, and results on this problem is discussed in Section~\ref{section5}.

After five years from the first announcement of the cohomological rigidity problem, many interesting results have been found, but still the complete solution to the problem is far from being settled. For this, some new insight might be needed, but we hope some enlightening new discoveries are made in a near future.

\section{Cohomological rigidity problem for toric manifolds}\label{section1}

A toric variety of complex dimension $n$ is a normal algebraic variety over the complex numbers $\C$ with an effective algebraic action of $(\C^*)^n$ having an open dense orbit, where $\C^*=\C\backslash\{0\}$. On the other hand, a fan of real dimension $n$ is a collection of cones in $\R^n$ with the origin as the vertex satisfying certain conditions.  A fundamental theorem in toric geometry says that there is a bijective correspondence between toric varieties of complex dimension $n$ and fans of real dimension $n$, see \cite{fult93} and \cite{oda88} for details.  Therefore, the classification of toric varieties \emph{as varieties} reduces to the classification of fans, and through this correspondence compact non-singular toric varieties, which we call \emph{toric manifolds} in this paper, are classified in some cases, see \cite[Section 3]{ma-su08}.

However, little is known about the diffeomorphism or homeomorphism classification of toric manifolds.  As is well-known, the cohomology ring $H^*(X;\Z)$ of a toric manifold $X$ is explicitly described in terms of the fan associated with $X$, i.e., $H^*(X;\Z)$ is the face ring of the underlying simplicial complex of the fan modulo a linear system of parameters, in particular, $H^*(X;\Z)$ is generated by degree two elements as a graded ring. See \cite[p.106]{fult93} or \cite[p.134]{oda88} for the precise description of the cohomology ring of toric manifolds. Cohomology ring does not distinguish closed smooth manifolds in general but we have a feeling that most of closed smooth manifolds do not have symmetry of a large torus.  Based on this feeling and by checking some examples, the second and third authors addressed the following problem in \cite{ma-su08}.

\begin{CRPT}
Are two toric manifolds diffeomorphic or homeomorphic if their cohomology rings with $\Z$-coefficients are isomorphic as graded rings?
\end{CRPT}

A toric manifold of complex dimension one is $\CP^1$.  Toric manifolds of complex dimension two are known to be $\CP^2$, Hirzebruch surfaces and their blow-ups (see \cite{fult93} or \cite{oda88}).  It is not difficult to check that the cohomological rigidity problem is affirmative for those manifolds.  The problem is open in higher dimensions but some partial affirmative solutions are obtained, which we shall survey below.

\subsection{Bott manifolds}
A \emph{Bott tower} of height $n$ is a sequence of $\CP^{1}$-bundles
\begin{equation} \label{eqn:Bott tower}
B_n\stackrel{\pi_n}\longrightarrow B_{n-1} \stackrel{\pi_{n-1}}\longrightarrow
\dots \stackrel{\pi_2}\longrightarrow B_1 \stackrel{\pi_1}\longrightarrow
B_0=\{\text{a point}\},
\end{equation}
where each $\pi_i\colon B_i\to B_{i-1}$ for $i=1,\dots,n$ is the projectivization of a Whitney sum of two complex line bundles over $B_{i-1}$. We call $B_i$ an \emph{$i$-stage Bott manifold} or a \emph{Bott manifold of height $i$}. It is not difficult to see that all Bott manifolds are toric manifolds. We are concerned with the diffeomorphism classification of $n$-stage Bott manifolds $B_n$.  Even if two Bott towers of height $n$ are different, their $n$-stage Bott manifolds $B_n$ can be diffeomorphic. For example, two-stage Bott manifolds $B_2$ are Hirzebruch surfaces, and there are only two diffeomorphism types, $(\CP^1)^2$ and $\C P^2\#\overline{\CP^2}$, among them although there are infinite choices of fibrations $\pi_2\colon B_2\to B_1=\CP^1$ for each diffeomorphism type of $B_2$.  The following theorem is proved in \cite{ch-ma-su10} when $n=3$ and recently proved by the first author in \cite{ch:preprint} when $n=4$.

\begin{theorem}[\cite{ch-ma-su10}, \cite{ch:preprint}]
The cohomological rigidity problem with diffeomorphism is affirmative for $n$-stage Bott manifolds when $n\le 4$.
\end{theorem}

If all the fibrations in \eqref{eqn:Bott tower} are trivial, then $B_n$ is diffeomorphic to $(\CP^1)^n$, and, hence, $H^*(B_n;\Z)\cong H^*((\CP^1)^n;\Z)$ as graded rings.  The following theorem shows that the converse holds.

\begin{theorem}[\cite{ma-pa08}] \label{theo:ma-pa08}
All the fibrations in \eqref{eqn:Bott tower} are trivial, in particular, $B_n$ is diffeomorphic to $(\CP^1)^n$ if $H^*(B_n;\Z)\cong H^*((\CP^1)^n;\Z)$ as graded rings.
\end{theorem}

We say that $B_n$ is \emph{$\Q$-trivial} if $H^*(B_n;\Q)\cong H^*((\CP^1)^n;\Q)$ as graded rings.  As remarked above, there are only two diffeomorphism types among Hirzebruch surfaces: one is $(\CP^1)^2$ and the other is $\CP^2\#\overline{\CP^2}$, and both of them can be easily seen to be $\Q$-trivial. This motivated us to study $\Q$-trivial Bott manifolds although there are infinitely many $n$-stage Bott manifolds which are not $\Q$-trivial when $n\ge 3$.

\begin{theorem}[\cite{ch-ma-pre}]
The cohomological rigidity problem with diffeomorphism is affirmative for $\Q$-trivial Bott manifolds.  Moreover, there are $p(n)$ number of diffeomorphism classes in $\Q$-trivial $n$-stage Bott manifolds, where $p(n)$ denotes the number of partitions of $n$.
\end{theorem}

We say that the Bott tower \eqref{eqn:Bott tower} is \emph{$k$-twisted} if there are exactly $k$ nontrivial topological fibrations in the tower.
The number $k$ is apparently an invariant of the tower. However, it is
not obvious that the number $k$ is independent of the choice of Bott
tower structures for a Bott manifold. Indeed, for a Bott manifold
$B_n$, there might be two different Bott towers whose last stage is
$B_n$. Therefore, the twist number $k$ for a Bott manifold may not be
well-defined. But, the twist number for a Bott manifold is, indeed,
well-defined, and is actually an invariant of $H^\ast(B_n)$, see \cite{ch-su-pre}.

\begin{theorem}[\cite{ch-su-pre}]
The cohomological rigidity problem with diffeomorphism is affirmative for Bott manifolds with one twist.
\end{theorem}

Since $\pi _i: B_i \to B_{i-1}$ has a cross section, the induced homomorphism $\pi _i^* : H^*(B_{i-1};\Z) \to H^*(B_i;\Z)$ is injective for each $i$. Therefore, $H^*(B_{i-1};\Z)$ can be regarded as a subring of $H^*(B_i;\Z)$ via $\pi _i^*$, so that we have a filtration by subrings:
\[
H^*(B_n;\Z)\supset H^*(B_{n-1};\Z)\supset \dots \supset H^*(B_1;\Z)\supset H^*(B_0;\Z).
\]
Let
\begin{equation} \label{another}
B'_n\stackrel{\pi'_n}\longrightarrow B'_{n-1} \stackrel{\pi'_{n-1}}\longrightarrow
\dots \stackrel{\pi'_2}\longrightarrow B'_1 \stackrel{\pi'_1}\longrightarrow
B'_0=\{\text{a point}\},
\end{equation}
be another Bott tower of height $n$.  We say that two Bott towers \eqref{eqn:Bott tower} and \eqref{another} are isomorphic if there is a collection of diffeomorphisms $\{f_i\colon B_i\to B'_i\}_{i=1}^n$ which commute with the projections $\{\pi_i\}_{i=1}^n$ and $\{\pi'_i\}_{i=1}^n$.  An isomorphism between Bott towers induces an isomorphism of their cohomology rings preserving the filtrations.  The following shows that the converse is also true.

\begin{theorem}[\cite{ishi-pre}]
Any graded ring isomorphism between the cohomology rings of $n$-stage Bott manifolds with $\Z$-coefficients preserving the filtrations can be realized by an isomorphism between the Bott towers of height $n$.
\end{theorem}

\subsection{Generalized Bott manifolds}

The notion of Bott tower and Bott manifold can naturally be generalized.  A \emph{generalized Bott tower} of height $n$ is a sequence of $\CP^{m_i}$-bundles with $m_i\ge 1$:
\begin{equation} \label{eqn:generalized Bott tower}
B_n\stackrel{\pi_n}\longrightarrow B_{n-1} \stackrel{\pi_{n-1}}\longrightarrow
\dots \stackrel{\pi_2}\longrightarrow B_1 \stackrel{\pi_1}\longrightarrow
B_0=\{\text{a point}\},
\end{equation}
where each $\pi_i\colon B_i\to B_{i-1}$ for $i=1,\dots,n$ is the projectivization of a Whitney sum of $m_i+1$ complex line bundles over $B_{i-1}$. We call $B_i$ an \emph{$i$-stage generalized Bott manifold} or a \emph{generalized Bott manifold of height $i$}.  They are also toric manifolds.  Similarly to the case of Bott manifolds, we are concerned with the diffeomorphism classification of $n$-stage generalized Bott manifolds $B_n$ and they can be diffeomorphic even if the towers are different.  The variety classification of two-stage generalized Bott manifolds is completed by Kleischmidt \cite{klei88} and their diffeomorphism classification is done by the authors in \cite{ch-ma-su10}.  As a corollary of the latter result we have

\begin{theorem}[\cite{ch-ma-su10}] \label{theo:2-stage}
The cohomological rigidity problem with diffeomorphism is affirmative for two-stage generalized Bott manifolds.
\end{theorem}

The paper \cite{ch-ma-su10} also generalizes Theorem~\ref{theo:ma-pa08} as follows.

\begin{theorem}[\cite{ch-ma-su10}] \label{thm:rigidity-ch-ma-su10}
All the fibrations in \eqref{eqn:generalized Bott tower} are trivial, in particular, $B_n$ is diffeomorphic to $\prod_{i=1}^n\CP^{m_i}$ if $H^*(B_n;\Z)\cong H^*(\prod_{i=1}^n\CP^{m_i};\Z)$ as graded rings.
\end{theorem}

Similarly to the case of Bott manifolds, we say that a generalized Bott manifold is \emph{$\Q$-trivial} if $H^*(B_n;\Q)\cong H^*(\prod_{i=1}^n\CP^{m_i};\Q)$ as graded rings.

\begin{theorem}[\cite{pa-su-pre}]
A $\Q$-trivial Bott manifold with $m_i \geq 2$ for all $i=1, \ldots, n$ is weakly equivariantly diffeomorphic to $\prod_{i=1}^n\CP^{m_i}$ with the standard torus action.
\end{theorem}

So far, the cohomological rigidity problem for toric manifolds seems not studied for other toric manifolds.

\section{Topological analogues of toric manifolds} \label{section2}

\subsection{Quasitoric manifolds}

Around 1990 Davis-Januszkiewicz \cite{da-ja91} introduced what is now called a \emph{quasitoric manifold}  (see \cite{bu-pa02}),  which is a topological analogue of a toric manifold.  A quasitoric manifold is a closed smooth manifold $M$ of even dimension, say $2n$, with an effective smooth action of $(S^1)^n$, satisfying the following two conditions:
\begin{enumerate}
\item $M$ is locally equivariantly diffeomorphic to a representation space of $(S^1)^n$,
\item the orbit space $M/(S^1)^n$ is a simple convex polytope.
\end{enumerate}
We say that $M$ is a quasitoric manifold over a simple polytope $Q$ when $Q=M/(S^1)^n$.

A toric manifold $X$ of complex dimension $n$ with the restricted action of $(S^1)^n$ is quasitoric in many cases.  In fact, this is the case when $X$ is projective and it is not known whether there is a toric manifold $X$ which is not quasitoric.  One can see that the restricted action of $(S^1)^n$ on $X$ satisfies condition (1) above so that the orbit space $X/(S^1)^n$ is a  manifold with corners.  Moreover, since the action of $(S^1)^n$ on $X$ is the restriction of the action of $(\C^*)^n$, the orbit space $X/(S^1)^n$ has the residual action of $(\C^*)^n/(S^1)^n\cong \R^n$.  This implies that all faces of $X/(S^1)^n$ (even $X/(S^1)^n$ itself) are contractible, and any intersection of faces is connected unless it is empty, so the orbit space looks like a simple polytope.  However, the orbit space might not be a simple polytope and there might exist a toric manifold which is not quasitoric.

On the other hand, there are many quasitoric manifolds which do not arise from toric manifolds.  For example, $\C P^2\#\C P^2$ with some smooth action of $(S^1)^2$ is quasitoric but not toric because it does not allow a complex (even almost complex) structure, as is well-known.

The cohomological rigidity problem is also asked for quasitoric manifolds in \cite{ma-su08}.  We note that an $n$-stage Bott manifold with the restricted compact torus action is a quasitoric manifold over an $n$-cube, more generally, an $n$-stage generalized Bott manifold in \eqref{eqn:generalized Bott tower} with the restricted compact torus action is a quasitoric manifold over $\prod_{i=1}^n\Delta^{m_i}$ a product of simplices where $\Delta^{m_i}$ denotes an $m_i$-simplex.  Most of quasitoric manifolds over a product of simplices are generalized Bott manifolds but there are quasitoric manifolds over a product of simplices which are not generalized Bott manifolds.  For instance, $\CP^2\#\CP^2$ with some smooth action of $(S^1)^2$ is a quasitoric manifold over a square but not toric as remarked before, in particular, not a Bott manifold.

\begin{theorem}[\cite{ch-pa-su-pre}]
Quasitoric manifolds over a product of two simplices are homeomorphic if their cohomology rings with $\Z$-coefficients are isomorphic. Moreover, if we fix a product of two simplices, denoted $Q$, then there are only finitely many quasitoric manifolds over $Q$ which are not generalized Bott manifolds.
\end{theorem}

\begin{remark}
The first statement in the theorem above depends on Proposition 1.8 in \cite{da-ja91} which holds up to equivariant homeomorphism. This is the reason why the first statement is stated as \emph{homeomorphic}.
\end{remark}

\subsection{Topological toric manifolds}

Recently, H.~Ishida, Y.~Fukukawa and the second author \cite{is-fu-ma-pre} have introduced another topological analogue of toric manifolds, called \emph{topological toric manifolds}.  It is known that a toric manifold of complex dimension $n$ is covered by finitely many invariant open subsets each equivariantly and algebraically isomorphic to a direct sum of complex one-dimensional \emph{algebraic} representation spaces of $(\C^*)^n$.  Based on this observation a topological toric manifold is defined as follows.

\begin{definition}
A closed smooth manifold $X$ of dimension $2n$ with an effective \emph{smooth} action of $(\C^*)^n$ having an open dense orbit is a (compact) \emph{topological toric manifold} if it is covered by finitely many invariant open subsets each equivariantly diffeomorphic to a direct sum of complex one-dimensional \emph{smooth} representation spaces of $(\C^*)^n$.
\end{definition}

To a topological toric manifold $X$, one can associate a combinatorial object $\Delta(X)$ called a complete non-singular topological fan similarly to the toric case.  The main theorem in \cite{is-fu-ma-pre} says that the correspondence $X\to \Delta(X)$ is a bijection between topological toric manifolds and complete non-singular topological fans.

Apparently from the definition, a toric manifold is a topological toric manifold.  But topological toric manifolds turn out to be much more abundant than toric manifolds because there are many more smooth representations of $(\C^*)^n$ than algebraic ones.  This stems from the fact that since $\C^*=\R_{>0}\times S^1$ as smooth groups, any \emph{smooth} endomorphism of $\C^*$ is of the form
\begin{equation} \label{intro1}
\text{$g\mapsto |g|^{b+\sqrt{-1}c}\big(\frac{g}{|g|}\big)^v$\quad with $(b+\sqrt{-1}c,v)\in \C\times\Z$}
\end{equation}
and this endomorphism is algebraic if and only if $b=v$ and $c=0$.  Therefore the group $\Hom(\C^*,\C^*)$ of \emph{smooth} endomorphisms of $\C^*$ is isomorphic to $\C\times\Z$ while the group $\Hom_{alg}(\C^*,\C^*)$ of \emph{algebraic} endomorphisms of $\C^*$ is isomorphic to $\Z$.

Nevertheless, topological toric manifolds have similar topological properties to toric manifolds.  For instance, the orbit space of a topological toric manifold $X$ by the restricted compact torus action is a manifold with corners whose faces (even the orbit space itself) are all contractible and any intersection of faces is connected unless it is empty, so the orbit space looks like a simple polytope similarly to the toric case.  However, it is shown in \cite{is-fu-ma-pre} that the orbit space is not necessarily a simple polytope and the family of topological toric manifolds with restricted compact torus actions is strictly larger than the family of quasitoric manifolds up to equivariant homeomorphism.  One can also see that  $H^*(X;\Z)$ has the same presentation as the toric or quasitoric case, i.e., a Stanley-Reisner ring modulo a linear system of parameters,  in particular, $H^*(X;\Z)$ is generated by degree two elements as a ring like the toric or quasitoric case.  Therefore, it is reasonable to ask the cohomological rigidity problem for topological toric manifolds.

\begin{CRPTT}
Are two topological toric manifolds diffeomorphic or homeomorphic if their cohomology rings with $\Z$-coefficients are isomorphic as graded rings?
\end{CRPTT}

So far, no counterexample to this problem is known.

\subsection{Torus manifolds}
Hattori and Masuda \cite{ha-ma03} introduced a \emph{torus manifold} (or \emph{unitary toric manifold} in an earlier terminology \cite{masu99}),  which is a closed, smooth manifold of dimension $2n$ with an effective smooth action of $(S^1)^n$ having a fixed point, and developed an analogous theory to toric geometry for torus manifolds to some extent.  Topological toric manifolds with restricted compact torus actions are torus manifolds, so the class of torus manifolds is most general.  However, this class seems too general to study.  For instance, there are many torus manifolds which are not simply connected and have non-trivial odd degree cohomology although topological toric manifolds are simply connected and their cohomology rings are generated by degree two elements as mentioned before.  Moreover, torus manifolds are not cohomologically rigid even if they are simply connected and have vanishing odd degree cohomology, as is shown in the following example.

\begin{example}[S.~Kuroki] \label{kuroki}
Let $M=\C P^1\times S^{2\ell}$ and $M'=S^3\times_{S^1} S(\C\oplus \R^{2\ell-1})$ where $S^1$ acts on $S^3$ freely in a natural way, on $\C$ as scalar multiplication, on $\R^{2\ell-1}$ trivially and $S(\cdot)$ stands for the unit sphere.  One can see easily that both $M$ and $M'$ admit smooth actions of $(S^1)^{\ell+1}$ so that they become torus manifolds.  Moreover, they are simply connected, have isomorphic cohomology rings with $\Z$-coefficients and their odd degree cohomology groups vanish.  However, $M$ is spin while $M'$ is non-spin as is easily checked, so they are not diffeomorphic, even not homotopy equivalent.
\end{example}

Nevertheless, the topological classification problem for torus manifolds is still interesting. Let $G$ be a compact, connected Lie group with maximal torus $(S^1)^n$.
It is shown in \cite{K2} that if a torus manifold has a codimension one extended $G$-action, then it is a product of even-dimensional spheres and the following manifolds
\begin{equation} \label{2types}
S^{2m +1} \times_{S^{1}} P(\C_\rho^k \oplus \C^\ell),\quad S^{2m +1} \times_{S^{1}} S(\C_\rho^{k} \oplus \R^{2\ell-1}) \quad (m,k, \ell \in \N, \rho \in \Z),
\end{equation}
where $P(\cdot)$ stands for projectivization, $S(\cdot)$ stands for the unit sphere as before, and $S^{1}$ acts on $S^{2m +1}$ naturally, and on $\C_{\rho}$ through the group homomorphism $S^{1}\to S^{1}$ defined by $t \mapsto t^\rho$.
We note that manifolds of the former type in \eqref{2types} are special two-stage generalized Bott manifolds, so they are distinguished by their cohomology rings up to diffeomorphism by Theorem~\ref{theo:2-stage}. Manifolds of the latter type in \eqref{2types} are two-stage generalized Bott manifolds if and only if $k=\ell=1$, and unless $k=\ell=1$, their cohomology rings are not generated by degree two elements.  Those manifolds are not necessarily distinguished by their cohomology rings up to diffeomorphism as is seen in Example~\ref{kuroki} and more invariants are necessary to distinguish them.

\begin{theorem}[\cite{ch-ku-pre}]
Two manifolds of the latter type in \eqref{2types} are diffeomorphic if there is a cohomology ring isomorphism preserving their Stiefel-Whitney classes and Pontrjagin classes.
\end{theorem}

Looking at this result, it is natural to ask the following.

\begin{problem}
What is the most general class in torus manifolds whose topological types are classified by their cohomology rings, their Stiefel-Whitney and Pontrjagin characteristic classes?
\end{problem}

\subsection{Symplectic toric manifolds}
A \emph{symplectic toric manifold} (or a \emph{toric symplectic manifold}) is a closed connected symplectic manifold of even dimension, say $2n$, with an effective Hamiltonian action of $(S^1)^n$.  A symplectic toric manifold is known to be equivariantly diffeomorphic to a projective toric manifold with restricted compact torus action (\cite{delz88}).  Recently, an analogous cohomological rigidity problem has been studied by D.~McDuff for symplectic toric manifolds.
\begin{theorem}[Borisov-McDuff \cite{mcdf10}]
Let $R$ be a commutative ring of finite rank with even grading, and write $R_{\R} := R\otimes_\Z\R$.  Suppose given elements $[\omega]\in R_{\R}$ and $c_1, c_2\in R$ of degrees 2, 2 and 4 respectively. Then, up to equivariant symplectomorphism, there are at most finitely many symplectic toric manifolds $(M,\omega)$ for which there is a ring isomorphism $\Psi\colon H^*(M;\Z)\to R$ that takes the symplectic class and the Chern classes $c_i(M)$, $i = 1, 2$, to the given elements $[\omega] \in R_{\R}, c_i\in R$.
\end{theorem}

\section{Real version}\label{section3}

\subsection{Real toric manifolds}

The set of real points in a toric manifold $X$ is called a real toric manifold and denoted by $X(\R)$.  For instance, when $X=\CP^n$, $X(\R)=\R P^n$. We may ask the cohomological rigidity problem for real toric manifolds, but $H^*(X(\R);\Z)$ seems complicated and unknown in general, so it seems not practical to ask the problem although the authors do not know any counterexample to the problem.  However, it is known that
\[
H^*(X(\R);\Z/2)\cong H^{2*}(X;\Z)\otimes_{\Z} \Z/2=H^{2*}(X;\Z/2),
\]
so that $H^*(X(\R);\Z/2)$ can be described in terms of the fan associated with $X$ since so is $H^*(X;\Z)$. Because of this reason the cohomological rigidity problem for real toric manifolds was addressed with $\Z/2$-coefficient in \cite{ka-ma09}, and it is shown there that the problem is affirmative for real Bott manifolds, where a \emph{real Bott manifold} is the real analogue of a Bott manifold obtained by replacing $\C$ in \eqref{eqn:Bott tower} by $\R$.

\begin{theorem}[\cite{ka-ma09}] \label{ka-ma}
Real Bott manifolds are diffeomorphic if their cohomology rings with $\Z/2$-coefficients are isomorphic as graded rings.
\end{theorem}

Real Bott manifolds provide examples of flat riemannian manifolds.  The theorem above has been reproved and improved in \cite{ch-ma-ou-pre}.  The paper \cite{ch-ma-ou-pre} also relates real Bott manifolds to acyclic digraphs (directed graphs with no directed cycle), and shows that the classification of real Bott manifolds of dimension $n$ up to diffeomorphism corresponds to the classification of non-isomorphic acyclic digraphs on $n$ vertices up to two graph operations where one is a well-known operation called a local complementation and the other seems a new operation.  Using this correspondence, the number of diffeomorphism classes of real Bott manifolds is explicitly counted up to dimension eight,  see Table~\ref{tab:Bott equivalence}.  Some of real Bott manifolds are orientable and some of orientable ones admit a symplectic form.  There are criteria of whether a real Bott manifold is orientable (\cite{ka-ma09}) or symplectic (\cite{ishi10}), and using those criteria, one can also count the number of diffeomorphism classes of orientable or symplectic real Bott manifolds for low dimensions, see also Table~\ref{tab:Bott equivalence}.

\begin{table}
    \centering
    \begin{tabular}{c|c c c c c c c c c c }
        \hline
        $n$ & 1 & 2 & 3 & 4  & 5 & 6 & 7 & 8&9&10\\ \hline
        $\Diff{n} $ & 1 & 2 & 4 & 12 & 54 & 472 & 8,512 & 328,416&?&?\\ \hline
        $\Ori{n}$ & 1 & 1 & 2 & 3 & 8 & 29 & 222 & 3,607&131,373&?\\ \hline
        $\Sym{n} $ &   & 1 &   & 2 &   & 6  &     &  31 &&416\\
        \hline
    \end{tabular}
    \caption{The numbers $\Diff{n},~\Ori{n},~\Sym{n}$ of $n$-dimensional real Bott manifolds, orientable real Bott manifolds and symplectic real Bott manifolds up to diffeomorphism, respectively.}
    \label{tab:Bott equivalence}
\end{table}

The paper \cite{ch-ma-ou-pre} also proves the following.

\begin{theorem}[Unique Decomposition Property] \label{UDP}
    The decomposition of a real Bott manifold into a product of indecomposable real Bott manifolds is unique up to permutations of the indecomposable factors.
\end{theorem}

Here a real Bott manifold is said to be indecomposable if it is not diffeomorphic to a product of two real Bott manifolds of positive dimension. Since $S^1$ is a real Bott manifold $\R P^1$, Theorem~\ref{UDP} implies the following corollary.

\begin{corollary}[Cancellation Property]
    Let $M$ and $M'$ be real Bott manifolds. If $S^1\times M$ and $S^1\times M'$ are diffeomorphic, then so are $M$ and $M'$.
\end{corollary}

As mentioned before, real Bott manifolds are compact flat riemannian manifolds.  We note that the cancellation property above does not hold for general compact flat riemannian manifolds (\cite{char65-1}).  These results imply that toric or real toric manifolds would have some rigidity properties, and it would be interesting to find them.

The real analogue of a generalized Bott manifold obtained by replacing $\C$ by $\R$ in \eqref{eqn:generalized Bott tower} is called a \emph{generalized real Bott manifold}.  Looking at Theorem~\ref{ka-ma} and Theorem~\ref{theo:2-stage}, the readers might expect that Theorem~\ref{ka-ma} can be extended to generalized real Bott manifolds, but the diffeomorphism classes of two-stage generalized real Bott manifolds are not necessarily distinguished by their cohomology rings with $\Z/2$-coefficients (\cite{masu10}).  However, if two two-stage generalized real Bott manifolds are homotopy equivalent, then they are diffeomorphic (\cite{masu10}).  Therefore, we may ask the following.

\begin{HRPRT} \label{HRPRT}
Are two real toric manifolds diffeomorphic or homeomorphic if they are homotopy equivalent?
\end{HRPRT}

Real toric manifolds are not simply connected and they are often aspherical manifolds.  Therefore, we may think of the problem above as a real toric version of the famous Borel conjecture on aspherical manifolds.

\subsection{Small covers}

The notion of small cover was introduced by Davis-Januszkiewicz (\cite{da-ja91}) as a real analogue of a quasitoric manifold.  It is a closed smooth manifold $N$ of dimension $n$ with an effective smooth action of $(S^0)^n$ satisfying the following two conditions:
\begin{enumerate}
\item $N$ is locally equivariantly diffeomorphic to a representation space of $(S^0)^n$,
\item the orbit space $N/(S^0)^n$ is a simple convex polytope.
\end{enumerate}
We say that $N$ is a small cover over a simple polytope $Q$ when $Q=N/(S^0)^n$.  Generalized real Bott manifolds are small covers over a product of simplices, and conversely, a small cover over a product of simplices is homeomorphic to a generalized real Bott manifold (\cite{ch-ma-su10}).  Similarly to the relation between toric manifolds and quasitoric manifolds, it is uncertain that any real toric manifold is a small cover but there are many small covers which are not real toric.

One can also ask the cohomological rigidity problem for small covers $N$.  Similarly to the real toric case, $H^*(N;\Z)$ is not well understood but $H^*(N;\Z/2)$ is explicitly described.  Therefore, the cohomological rigidity problem is asked for small covers with $\Z/2$-coefficient as before. This problem is affirmative for small covers over an $n$-cube (i.e., real Bott manifolds) by Theorem~\ref{ka-ma} but is not always affirmative because two-stage generalized real Bott manifolds are not distinguished by their cohomology rings with $\Z/2$-coefficients as mentioned before.  However, the following is known.

\begin{theorem}[\cite{ca-lu-pre}]
Small covers over a product of an $m$-gon and an interval are distinguished up to homeomorphism by their cohomology rings with $\Z/2$-coefficients.
\end{theorem}

\begin{remark}
We may replace \lq\lq homeomorphism" above by \lq\lq diffeomorphism" since they are equivalent in dimension three as is well known.
\end{remark}

Looking at these results, it is natural to ask the following.

\begin{problem}
For which simple polytope $Q$ are small covers over $Q$ distinguished up to homeomorphism or diffeomorphism by their cohomology rings with $\Z/2$-coefficients?
\end{problem}

Finally, we remark that it is not known whether two small covers are homeomorphic or diffeomorphic if they
are homotopy equivalent.

\subsection{Real topological toric manifolds}

Real analogue of topological toric manifolds has also been introduced in \cite{is-fu-ma-pre}.

\begin{definition}
A closed smooth manifold $Y$ of dimension $n$ with an effective \emph{smooth} action of $(\R^*)^n$ having an open dense orbit is called a (compact) \emph{real topological toric manifold} if it is covered by finitely many invariant open subsets each equivariantly diffeomorphic to a direct sum of real one-dimensional \emph{smooth} representation spaces of $(\R^*)^n$.
\end{definition}

Similarly to the complex case, the family of real topological toric manifolds properly contains the family of real toric manifolds.  Moreover, the family of real topological toric manifolds with the restricted action of  the 2-torus group $(S^0)^n$ is strictly larger than the family of small covers up to equivariant homeomorphism.  Therefore, it is reasonable to ask the rigidity problems discussed above for real topological toric manifolds.

\section{Strong cohomological rigidity problems}\label{section4}

The cohomological rigidity problem for toric manifolds asks whether two toric manifolds are diffeomorphic or homeomorphic if their cohomology rings are isomorphic. More strongly, we may ask the following (\cite{ma-su08}).

\begin{SCRPT}
Let $M$ and $N$ be toric manifolds, and $\varphi \colon H^\ast(N;\Z) \to H^\ast(M;\Z)$ a ring isomorphism. Then, does there exist a diffeomorphism or a homeomorphism $f \colon M \to N$ such that $f^\ast = \varphi$?
\end{SCRPT}

One can also ask the problem for quasitoric manifolds or topological toric manifolds.
In general, the answer to this problem with diffeomorphism is ``no''.  In fact, Friedman and Morgan \cite{Fr-Mo-1988} show that some cohomology ring automorphism of $\CP^2 \# 10 \overline{\CP^2}$ is not induced by its diffeomorphism while it is a (projective) toric manifold.  
To the contrary, we have some partial affirmative solutions to the problem with diffeomorphism as follows.

\begin{theorem}
The strong cohomological rigidity problem with diffeomorphism is affirmative in the following cases:
\begin{enumerate}
\item either $M$ or $N$ is $\prod_{i=1}^m \CP^{n_i}$ {\rm (\cite{ma-pa08} when $n_i=1$ for any $i$, \cite{ch-su-pre2} in general)},
\item both $M$ and $N$ are  $\Q$-trivial Bott manifolds {\rm (\cite{ch-ma-pre})},
\item both $M$ and $N$ are  three-stage Bott manifolds {\rm (\cite{ch:preprint})}.
\end{enumerate}
\end{theorem}

Note that Stiefel-Whitney classes are homotopy invariants, and Pontrjagin classes are homeomorphism invariants for closed smooth manifolds with torsion-free cohomology rings.
\begin{proposition}[\cite{ch-ma-su10}]
Let $M, M'$ be connected closed manifolds of dimension $n$ having isomorphic cohomology rings with $\Z/2$-coefficients. Suppose that $H^*(M; \Z/2)$ is generated by $H^r(M; \Z/2)$ for some $r$ as a ring. Then any ring isomorphism between $M$ and $M'$ preserves their Stiefel-Whitney classes.
\end{proposition}
Recently, the first author proved the following.
\begin{theorem}[\cite{ch:arxiv}]
Any graded ring isomorphism between the integral cohomology rings of $n$-stage Bott manifolds preserves their Pontrjagin classes. In particular, there are only finitely many diffeomorphism classes of Bott manifolds which are homotopy equivalent.
\end{theorem}


The strong cohomological rigidity problem can also be asked for real toric manifolds, small covers or real topological toric manifolds.  We may consider the problem for those manifolds with  $\Z$-coefficient as above but it would be practical to take $\Z/2$-coefficient as mentioned before.  Since the cohomological rigidity problem with $\Z/2$-coefficient is negative for generalized real Bott manifolds as we stated before, the strong cohomological rigidity with $\Z/2$-coefficient fails for them.  However, it holds for real Bott manifolds.
\begin{theorem}[\cite{ch-ma-ou-pre}]
Any isomorphism between the cohomology rings of two real Bott manifolds with $\Z/2$-coefficient is induced by a diffeomorphism.
\end{theorem}


\section{Rigidity problems for polytopes}\label{section5}
As mentioned before, a complete non-singular fan is associated with a toric manifold. One of interesting questions is to ask whether the cohomology ring of a toric manifold $X$ determines the combinatorial structure of the underlying simplicial complex $\Sigma_X$ of the fan associated with $X$. Although $H^\ast(X;\Z)$ does not determine $\Sigma_X$ in general, it sometimes does; for instance, if $H^\ast(X;\Z)$ is isomorphic to the cohomology ring of some Bott manifold, then $\Sigma_X$ is the boundary complex of a crosspolytope (\cite{ma-pa08}). Motivated by this, several notions of rigidity for simplicial complexes or polytopes have been introduced.  Here are two of them.
\begin{enumerate}
\item A simplicial complex $K$ is \emph{(cohomologically) rigid}  if  there is a toric manifold $X$ with $K=\Sigma_X$, and whenever there exists a toric manifold $Y$ with a graded ring isomorphism $H^\ast(X;\Z) \cong H^\ast(Y;\Z)$, $\Sigma_Y$ is isomorphic to $K$ as simplicial complexes (\cite{ma-su08}).
\item A simple polytope $P$ is \emph{cohomologically rigid} if there exists a quasitoric manifold $M$ over $P$, and whenever there exists a quasitoric manifold $N$ over a simple polytope $Q$ with a graded ring isomorphism $H^*(M;\Z) \cong H^*(N;\Z)$, $Q$ is combinatorially equivalent to $P$ (\cite{ch-pa-su10}).  A simplicial polytope is said to be \emph{cohomologically rigid} if its dual simple polytope is cohomologically rigid.
\end{enumerate}

The above rigidity of simplicial complexes or polytopes requires the existence of supporting toric manifolds or quasitoric manifolds. It raises some ambiguity. For instance, 
the dodecahedron is cohomologically rigid in the sense of (2) above (\cite{ch-pa-su10}), but we do not know whether the boundary complex of its dual polytope is rigid in the sense of (1) above because we do not know whether it supports a toric manifold. Moreover, we may use topological toric manifolds instead of toric manifolds or quasitoric manifolds to define the rigidity above. But then, this new notion of rigidity might be different from the rigidity in (1) or (2) above.  Thus, we are led to define rigidity it a purely algebraic and combinatorial way.

\subsection{Combinatorial rigidity}
Let $\k$ be a commutative ring with unity and $A=\k[v_1,\ldots, v_m]$ the polynomial graded ring in $v_1,\ldots, v_m$ over $\k$ with $\deg v_i=2$ for all $i$. Let $K$ be a simplicial complex with $\{1,\dots,m\}$ as the set of vertices. The \emph{Stanley-Reisner ring} $\k(K)$ of $K$ over $\k$ is defined to be $A/I_K$, where $I_K$ is the ideal generated by the monomials $v_{i_1}v_{i_2}\cdots v_{i_r}$ such that $\{i_1,i_2,\ldots,i_r\}$ does not support a face of $K$.

It is known that the simplicial complex $\Sigma_X$ associated with a toric manifold $X$ is Cohen-Macaulay, and there exists a linear system of parameters in $\Z(\Sigma_X)$. As mentioned before, $H^\ast(X;\Z)$ is the face ring $\Z(\Sigma_X)$ modulo some linear system of parameters.

\begin{definition}
Let $\K$ be a set of Cohen-Macaulay complexes. An element $K$ of $\K$ is \emph{cohomologically rigid} in $\K$ if $K' \in \K$ is isomorphic to $K$ whenever $\Z(K')/J'$ is isomorphic to $\Z(K)/J$ as graded rings where $J'$ and $J$ are the ideals generated by linear systems of parameters in $\Z(K')$ and $\Z(K)$, respectively.
\end{definition}

Let $\K_{CM}$ be the set of Cohen-Macaulay complexes, $\K_T$ the set of simplicial complexes supporting toric manifolds, and $\K_{poly}$ the set of the boundary complexes of simplicial polytopes. We remark that both $\K_T$ and $\K_{poly}$ are subsets of $\K_{CM}$.

\begin{problem}
Find all simplicial complexes which are cohomologically rigid in $\K_{CM}$, $\K_T$ or $\K_{poly}$.
\end{problem}

On the other hand, the cohomological rigidity of polytopes can be considered as a purely combinatorial problem. Let $K$ be a Cohen-Macaulay complex.
A \emph{finite free resolution}$[F:\phi]$ of $\Q(K)$ is an exact sequence
\begin{equation}\label{eq:ffr}
\xymatrix{ 0 \ar[r]& F_r \ar[r]^{\phi_r}& F_{r-1} \ar[r]^{\phi_{r-1}}&
  \cdots \ar[r]^{\phi_{2}} & F_1\ar[r]^{\phi_1}& F_0 \ar[r]^{\phi_0}&
  \Q(K) \ar[r] & 0},
\end{equation}
where $F_i$ is a finite free $\Q[x_1,\ldots,x_n]$-module and each $\phi_i$ is degree-preserving. If we take $F_{i}$ to be the module generated by the minimal basis of $\Ker(\phi_{i-1})$, then we get a \emph{minimal resolution} of $\Q(K)$. Since $\Q(K)$ is graded, so are all $F_i$'s, that is, $F_i = \bigoplus_{j}F_{i,2j}$. We denote by
$$
    \beta_{i,2j}(K) = \dim_\Q F_{i,2j}
$$ and call it the $(i,2j)$-th \emph{bigraded Betti number} of $K$. Choi-Panov-Suh \cite{ch-pa-su10} showed that $\Z(K)/J \cong \Z(K')/J'$ implies that $\beta_{i,j}(K)=\beta_{i,j}(K')$ for all $i,j$, where $K' \in \K_{CM}$, and $J \subset \Z(K)$ and $J' \subset \Z(K')$ are ideals generated by linear systems of parameters. This theorem raises the following new notion of rigidity for simplicial complexes.

\begin{definition}
Let $\K$ be a subset of Cohen-Macaulay complexes. An element $K$ of $\K$ is \emph{combinatorially rigid} in $\K$ if $K' \in \K$ is isomorphic to $K$ whenever $\beta_{i,j}(K')=\beta_{i,j}(K)$ for all $i,j$.
\end{definition}

The word `combinatorial rigidity' was named by the first author and J.~S.~Kim \cite{ch-ki:arxiv}. Note that if $K$ is combinatorially rigid in $\K$, then $K$ is cohomologically rigid in $\K$.
\begin{problem} \label{problem:ComRP}
Find all simplicial complexes which are combinatorially rigid in $\K_{CM}$, $\K_T$ or $\K_{poly}$.
\end{problem}

We note that the bigraded Betti numbers can be computed in combinatorial way \cite{Ho}. 
Hence, Problem~\ref{problem:ComRP} is a purely combinatorial problem. From now on, we consider only the combinatorial rigidity in $\K_{poly}$. We sometimes do not distinguish between a simplicial polytope and its boundary complex.

\subsection{Irreducible rigid polytopes}
Let $K_1$ and $K_2$ be simplicial polytopes. A connected sum of $K_1$ and $K_2$ is a polytope obtained by attaching a facet of $K_1$ and a facet of $K_2$. It depends on the way of choosing the two facets and identifying their vertices. A simplicial polytope $K$ is called \emph{reducible} if it can be expressed as a connected sum of two polytopes, and is called \emph{irreducible} otherwise. Let $C(K_1 \sharp K_2)$ denote the set of connected sums of $K_1$ and $K_2$. If there is only one connected sum of $K_1$ and $K_2$ up to isomorphism, then we will write the unique polytope as $K_1 \sharp K_2$.

\begin{example}
The polytopes in Table~\ref{tab:irr_rigid} are both irreducible and combinatorially rigid in $\K_{poly}$. For reader's convenience, we represent the polytopes in terms of both simple and simplicial terminologies. Furthermore, the dual of triangle-free $n$-dimensional simple polytopes with at most $2n+2$ facets, and the dual of all irreducible three-dimensional simple polytopes with at most $9$ facets are both irreducible and combinatorially rigid (\cite{ch-pa-su10}).

\begin{table}
  \centering
\begin{tabular}{|c|c|c|c|}
  \hline
  Simple & Simplicial & dimension & reference \\ \hline
  $\prod_{i=1}^m \Delta^{n_i}$ & $\partial \Delta^{n_1} \ast \cdots \ast \partial \Delta^{n_m}$ & $n = n_1 + \cdots + n_m$ & \cite{ch-pa-su10} \\
  $k$-gonal prism & $k$-gonal bi-pyramid & $3$ & \cite{ch-ki:arxiv} \\
  $k$-gonal edge-cut-prism & $k$-gonal semi bi-pyramid & $3$ & \cite{ch-ki:arxiv} \\
  dodecahedron & icosahedron &  $3$ & \cite{ch-pa-su10} \\
  \hline
\end{tabular}
  \caption{Irreducible combinatorially rigid polytopes}\label{tab:irr_rigid}
\end{table}
\end{example}

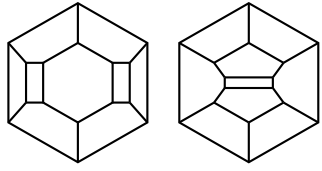
\begin{figure}
\psset{unit=10pt}
\begin{pspicture}(0,0)(6,6)
\pspolygon(0.40, 4.50)(0.40, 1.50)(3.00, 0.00)(5.60, 1.50)(5.60, 4.50)(3.00, 6.00)
\pspolygon(1.70, 3.75)(1.70, 2.25)(3.00, 1.50)(4.30, 2.25)(4.30, 3.75)(3.00, 4.50)
\psline(1.06, 3.75)(1.06, 2.25)\psline(4.94, 2.25)(4.94, 3.75)
\psline(1.06, 3.75)(1.70, 3.75)\psline(1.06, 3.75)(0.40, 4.50)
\psline(1.06, 2.25)(0.40, 1.50)\psline(1.06, 2.25)(1.70, 2.25)
\psline(4.94, 2.25)(5.60, 1.50)\psline(4.94, 2.25)(4.30, 2.25)
\psline(4.94, 3.75)(5.60, 4.50)\psline(4.94, 3.75)(4.30, 3.75)
\psline(3.00, 6.00)(3.00, 4.50)\psline(3.00, 0.00)(3.00, 1.50)
\end{pspicture}
\begin{pspicture}(0,0)(6,6)
\pspolygon(0.40, 4.50)(0.40, 1.50)(3.00, 0.00)(5.60, 1.50)(5.60, 4.50)(3.00, 6.00)
\psline(0.40, 4.50)(1.70, 3.75)\psline(0.40, 1.50)(1.70, 2.25)\psline(3.00, 0.00)(3.00, 1.50)
\psline(5.60, 1.50)(4.30, 2.25)\psline(5.60, 4.50)(4.30, 3.75)\psline(3.00, 6.00)(3.00, 4.50)
\psline(2.10, 2.80)(3.90, 2.80)\psline(2.10, 3.20)(3.90, 3.20)\psline(2.10, 2.80)(2.10, 3.20)
\psline(3.90, 2.80)(3.90, 3.20)\psline(2.10, 2.80)(1.70, 2.25)\psline(2.10, 3.20)(1.70, 3.75)
\psline(3.90, 2.80)(4.30, 2.25)\psline(3.90, 3.20)(4.30, 3.75)
\psline(1.70, 2.25)(3.00, 1.50)\psline(3.00, 1.50)(4.30, 2.25)\psline(4.30, 3.75)(3.00, 4.50)\psline(3.00, 4.50)(1.70, 3.75)
\end{pspicture}
  \caption{Examples of combinatorially non-rigid simple polytopes}\label{fig:non-rigid irreducible}
\end{figure}

On the other hand, not all irreducible polytopes are combinatorially rigid. It is known that there are $10$ irreducible three-dimensional simple polytopes with $10$ facets. Among them, the dual of two polytopes in Figure~\ref{fig:non-rigid irreducible} have the identical bigraded Betti numbers \cite{choi_y:2009}. So they are not combinatorially rigid in $\K_{poly}$. However, we do not know whether they are cohomologically rigid or not.

\begin{problem}
Find a polytope which is rigid cohomologically but not combinatorially in $\K_{poly}$.
\end{problem}

\subsection{Reducible rigid polytopes}
We have many examples of reducible non-rigid simplicial polytopes. A \emph{stacked polytope} is a simplicial polytope obtained by a sequence of connected sum of simplices from a simplex. The bigraded Betti numbers of a stacked polytope only depend on the number of simplices \cite{ch-ki:2010,te-hi:1997}. One can see easily that there are at least two elements in the set of connected sums of $k$ copies of $n$-simplices, where $n \geq 3$ and $k \geq 4$. Hence, they provide a lot of non-rigid simplicial polytopes combinatorially (even cohomologically, see \cite{ma-su08}) but reducible.

In many cases, the bigraded Betti numbers of the connected sum of simplicial polytopes is independent of the ways of connected sum, such as simplices, and two- and three-polytopes. In these cases, the rigidity of a reducible polytope $K \in C(K_1 \sharp K_2)$ is closely related to the singleton-ness of $C(K_1 \sharp K_2)$. Roughly speaking, reducible simplicial polytopes are combinatorially non-rigid unless both $K_1$ and $K_2$ admit sufficient symmetry. In two-dimensional case, a $k$-gon can be represented by the connected sum of $k-2$ triangles, and is combinatorially rigid. However, in higher dimension, there might be only a few examples which are both combinatorially rigid and reducible.

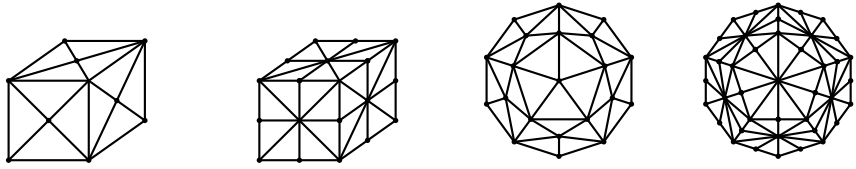
\begin{figure}
  \centering
  \psset{unit=30pt}
\begin{pspicture}(0,0)(3,1.5) \vput(0,1)0 \vput(0,0)0 \vput(1,0)0
  \vput(1.7,0.5)0 \vput(1.7,1.5)0 \vput(0.7,1.5)0
  \pspolygon(0,1)(0,0)(0,0)(1,0)(1,0)(1.7,0.5)(1.7,0.5)(1.7,1.5)(1.7,1.5)(0.7,1.5)
  \psline(1,1)(0,1) \psline(1,1)(1.7,1.5) \psline(1,1)(1,0)
  \vput(0.5,0.5)0 \psline(0.5,0.5)(1,1) \psline(0.5,0.5)(0,1)
  \psline(0.5,0.5)(0,0) \psline(0.5,0.5)(1,0) \vput(0.85,1.25)0
  \psline(0.85,1.25)(1,1) \psline(0.85,1.25)(0,1)
  \psline(0.85,1.25)(0.7,1.5) \psline(0.85,1.25)(1.7,1.5)
  \vput(1.35,0.75)0 \psline(1.35,0.75)(1,1) \psline(1.35,0.75)(1,0)
  \psline(1.35,0.75)(1.7,0.5) \psline(1.35,0.75)(1.7,1.5)
\end{pspicture}
\begin{pspicture}(0,0)(2,1.5) \vput(0,1)0 \vput(0,0)0 \vput(1,0)0
  \vput(1.7,0.5)0 \vput(1.7,1.5)0 \vput(0.7,1.5)0
  \pspolygon(0,1)(0,0)(0,0)(1,0)(1,0)(1.7,0.5)(1.7,0.5)(1.7,1.5)(1.7,1.5)(0.7,1.5)
  \psline(1,1)(0,1) \psline(1,1)(1.7,1.5) \psline(1,1)(1,0)
  \vput(0.5,0.5)0 \psline(0.5,0.5)(1,1) \psline(0.5,0.5)(0,1)
  \psline(0.5,0.5)(0,0) \psline(0.5,0.5)(1,0) \vput(0.5,1)0
  \psline(0.5,0.5)(0.5,1) \vput(0,0.5)0 \psline(0.5,0.5)(0,0.5)
  \vput(0.5,0)0 \psline(0.5,0.5)(0.5,0) \vput(1,0.5)0
  \psline(0.5,0.5)(1,0.5) \vput(0.85,1.25)0 \psline(0.85,1.25)(1,1)
  \psline(0.85,1.25)(0,1) \psline(0.85,1.25)(0.7,1.5)
  \psline(0.85,1.25)(1.7,1.5) \vput(0.5,1)0 \psline(0.85,1.25)(0.5,1)
  \vput(0.35,1.25)0 \psline(0.85,1.25)(0.35,1.25) \vput(1.2,1.5)0
  \psline(0.85,1.25)(1.2,1.5) \vput(1.35,1.25)0
  \psline(0.85,1.25)(1.35,1.25) \vput(1.35,0.75)0
  \psline(1.35,0.75)(1,1) \psline(1.35,0.75)(1,0)
  \psline(1.35,0.75)(1.7,0.5) \psline(1.35,0.75)(1.7,1.5)
  \vput(1,0.5)0 \psline(1.35,0.75)(1,0.5) \vput(1.35,0.25)0
  \psline(1.35,0.75)(1.35,0.25) \vput(1.7,1)0
  \psline(1.35,0.75)(1.7,1) \vput(1.35,1.25)0
  \psline(1.35,0.75)(1.35,1.25)
\end{pspicture}
\hspace{.5cm}
  \psset{unit=15pt}
\begin{pspicture}(-2,-2)(2,2) \vput(1.14127,0.37082)0 \vput(0,1.2)0
  \vput(-1.14127,0.37082)0 \vput(-0.705342,-0.97082)0
  \vput(0.705342,-0.97082)0
  \pspolygon(1.14127,0.37082)(0,1.2)(-1.14127,0.37082)(-0.705342,-0.97082)(0.705342,-0.97082)
  \vput(1.80701,0.587132)0 \vput(1.11679,1.53713)0 \vput(0,1.9)0
  \vput(-1.11679,1.53713)0 \vput(-1.80701,0.587132)0
  \vput(-1.80701,-0.587132)0 \vput(-1.11679,-1.53713)0 \vput(0,-1.9)0
  \vput(1.11679,-1.53713)0 \vput(1.80701,-0.587132)0
  \pspolygon(1.80701,0.587132)(1.11679,1.53713)(0,1.9)(-1.11679,1.53713)(-1.80701,0.587132)(-1.80701,-0.587132)(-1.11679,-1.53713)(0,-1.9)(1.11679,-1.53713)(1.80701,-0.587132)
  \psline(1.14127,0.37082)(1.80701,0.587132) \psline(0,1.2)(0,1.9)
  \psline(-1.14127,0.37082)(-1.80701,0.587132)
  \psline(-0.705342,-0.97082)(-1.11679,-1.53713)
  \psline(0.705342,-0.97082)(1.11679,-1.53713) \vput(0,0)0
  \psline(0,0)(1.14127,0.37082) \psline(0,0)(0,1.2)
  \psline(0,0)(-1.14127,0.37082) \psline(0,0)(-0.705342,-0.97082)
  \psline(0,0)(0.705342,-0.97082) \vput(0.822899,1.13262)0
  \psline(0.822899,1.13262)(0,1.2)
  \psline(0.822899,1.13262)(1.14127,0.37082)
  \psline(0.822899,1.13262)(1.80701,0.587132)
  \psline(0.822899,1.13262)(1.11679,1.53713)
  \psline(0.822899,1.13262)(0,1.9) \vput(-0.822899,1.13262)0
  \psline(-0.822899,1.13262)(-1.14127,0.37082)
  \psline(-0.822899,1.13262)(0,1.2) \psline(-0.822899,1.13262)(0,1.9)
  \psline(-0.822899,1.13262)(-1.11679,1.53713)
  \psline(-0.822899,1.13262)(-1.80701,0.587132)
  \vput(-1.33148,-0.432624)0
  \psline(-1.33148,-0.432624)(-0.705342,-0.97082)
  \psline(-1.33148,-0.432624)(-1.14127,0.37082)
  \psline(-1.33148,-0.432624)(-1.80701,0.587132)
  \psline(-1.33148,-0.432624)(-1.80701,-0.587132)
  \psline(-1.33148,-0.432624)(-1.11679,-1.53713) \vput(0,-1.4)0
  \psline(0,-1.4)(0.705342,-0.97082)
  \psline(0,-1.4)(-0.705342,-0.97082)
  \psline(0,-1.4)(-1.11679,-1.53713) \psline(0,-1.4)(0,-1.9)
  \psline(0,-1.4)(1.11679,-1.53713) \vput(1.33148,-0.432624)0
  \psline(1.33148,-0.432624)(1.14127,0.37082)
  \psline(1.33148,-0.432624)(0.705342,-0.97082)
  \psline(1.33148,-0.432624)(1.11679,-1.53713)
  \psline(1.33148,-0.432624)(1.80701,-0.587132)
  \psline(1.33148,-0.432624)(1.80701,0.587132)
\end{pspicture}
\hspace{.5cm}
\begin{pspicture}(-2,-2)(2,2) \vput(1.14127,0.37082)0 \vput(0,1.2)0
  \vput(-1.14127,0.37082)0 \vput(-0.705342,-0.97082)0
  \vput(0.705342,-0.97082)0
  \pspolygon(1.14127,0.37082)(0,1.2)(-1.14127,0.37082)(-0.705342,-0.97082)(0.705342,-0.97082)
  \vput(1.80701,0.587132)0 \vput(1.11679,1.53713)0 \vput(0,1.9)0
  \vput(-1.11679,1.53713)0 \vput(-1.80701,0.587132)0
  \vput(-1.80701,-0.587132)0 \vput(-1.11679,-1.53713)0 \vput(0,-1.9)0
  \vput(1.11679,-1.53713)0 \vput(1.80701,-0.587132)0
  \pspolygon(1.80701,0.587132)(1.11679,1.53713)(0,1.9)(-1.11679,1.53713)(-1.80701,0.587132)(-1.80701,-0.587132)(-1.11679,-1.53713)(0,-1.9)(1.11679,-1.53713)(1.80701,-0.587132)
  \psline(1.14127,0.37082)(1.80701,0.587132) \psline(0,1.2)(0,1.9)
  \psline(-1.14127,0.37082)(-1.80701,0.587132)
  \psline(-0.705342,-0.97082)(-1.11679,-1.53713)
  \psline(0.705342,-0.97082)(1.11679,-1.53713) \vput(0,0)0
  \psline(0,0)(1.14127,0.37082) \psline(0,0)(0,1.2)
  \psline(0,0)(-1.14127,0.37082) \psline(0,0)(-0.705342,-0.97082)
  \psline(0,0)(0.705342,-0.97082) \vput(0.570634,0.78541)0
  \psline(0,0)(0.570634,0.78541) \vput(-0.570634,0.78541)0
  \psline(0,0)(-0.570634,0.78541) \vput(-0.923305,-0.3)0
  \psline(0,0)(-0.923305,-0.3) \vput(0,-0.97082)0
  \psline(0,0)(0,-0.97082) \vput(0.923305,-0.3)0
  \psline(0,0)(0.923305,-0.3) \vput(0.822899,1.13262)0
  \psline(0.822899,1.13262)(0,1.2)
  \psline(0.822899,1.13262)(1.14127,0.37082)
  \psline(0.822899,1.13262)(1.80701,0.587132)
  \psline(0.822899,1.13262)(1.11679,1.53713)
  \psline(0.822899,1.13262)(0,1.9) \vput(0.570634,0.78541)0
  \psline(0.822899,1.13262)(0.570634,0.78541) \vput(1.47414,0.478976)0
  \psline(0.822899,1.13262)(1.47414,0.478976) \vput(1.4619,1.06213)0
  \psline(0.822899,1.13262)(1.4619,1.06213) \vput(0.558396,1.71857)0
  \psline(0.822899,1.13262)(0.558396,1.71857) \vput(0,1.55)0
  \psline(0.822899,1.13262)(0,1.55) \vput(-0.822899,1.13262)0
  \psline(-0.822899,1.13262)(-1.14127,0.37082)
  \psline(-0.822899,1.13262)(0,1.2) \psline(-0.822899,1.13262)(0,1.9)
  \psline(-0.822899,1.13262)(-1.11679,1.53713)
  \psline(-0.822899,1.13262)(-1.80701,0.587132)
  \vput(-0.570634,0.78541)0
  \psline(-0.822899,1.13262)(-0.570634,0.78541) \vput(0,1.55)0
  \psline(-0.822899,1.13262)(0,1.55) \vput(-0.558396,1.71857)0
  \psline(-0.822899,1.13262)(-0.558396,1.71857)
  \vput(-1.4619,1.06213)0 \psline(-0.822899,1.13262)(-1.4619,1.06213)
  \vput(-1.47414,0.478976)0
  \psline(-0.822899,1.13262)(-1.47414,0.478976)
  \vput(-1.33148,-0.432624)0
  \psline(-1.33148,-0.432624)(-0.705342,-0.97082)
  \psline(-1.33148,-0.432624)(-1.14127,0.37082)
  \psline(-1.33148,-0.432624)(-1.80701,0.587132)
  \psline(-1.33148,-0.432624)(-1.80701,-0.587132)
  \psline(-1.33148,-0.432624)(-1.11679,-1.53713)
  \vput(-0.923305,-0.3)0 \psline(-1.33148,-0.432624)(-0.923305,-0.3)
  \vput(-1.47414,0.478976)0
  \psline(-1.33148,-0.432624)(-1.47414,0.478976) \vput(-1.80701,0)0
  \psline(-1.33148,-0.432624)(-1.80701,0) \vput(-1.4619,-1.06213)0
  \psline(-1.33148,-0.432624)(-1.4619,-1.06213)
  \vput(-0.911067,-1.25398)0
  \psline(-1.33148,-0.432624)(-0.911067,-1.25398) \vput(0,-1.4)0
  \psline(0,-1.4)(0.705342,-0.97082)
  \psline(0,-1.4)(-0.705342,-0.97082)
  \psline(0,-1.4)(-1.11679,-1.53713) \psline(0,-1.4)(0,-1.9)
  \psline(0,-1.4)(1.11679,-1.53713) \vput(0,-0.97082)0
  \psline(0,-1.4)(0,-0.97082) \vput(-0.911067,-1.25398)0
  \psline(0,-1.4)(-0.911067,-1.25398) \vput(-0.558396,-1.71857)0
  \psline(0,-1.4)(-0.558396,-1.71857) \vput(0.558396,-1.71857)0
  \psline(0,-1.4)(0.558396,-1.71857) \vput(0.911067,-1.25398)0
  \psline(0,-1.4)(0.911067,-1.25398) \vput(1.33148,-0.432624)0
  \psline(1.33148,-0.432624)(1.14127,0.37082)
  \psline(1.33148,-0.432624)(0.705342,-0.97082)
  \psline(1.33148,-0.432624)(1.11679,-1.53713)
  \psline(1.33148,-0.432624)(1.80701,-0.587132)
  \psline(1.33148,-0.432624)(1.80701,0.587132) \vput(0.923305,-0.3)0
  \psline(1.33148,-0.432624)(0.923305,-0.3) \vput(0.911067,-1.25398)0
  \psline(1.33148,-0.432624)(0.911067,-1.25398)
  \vput(1.4619,-1.06213)0 \psline(1.33148,-0.432624)(1.4619,-1.06213)
  \vput(1.80701,0)0 \psline(1.33148,-0.432624)(1.80701,0)
  \vput(1.47414,0.478976)0
  \psline(1.33148,-0.432624)(1.47414,0.478976)
\end{pspicture}
 \caption{$\xi_1(C_8)$, $\xi_2(C_8)$, $\xi_1(D_{20})$ and $\xi_2(D_{20})$.}
  \label{fig:dodeca}
\end{figure}

The first author and J.S.Kim \cite{ch-ki:arxiv} studied the rigidity of reducible polytope of dimension three. They established a necessary condition for $C(K_1 \sharp K_2)$ to be a singleton, and find some reducible rigid polytopes of dimension three.
Let $T_4$, $C_8$, $O_6$, $D_{20}$ and $I_{12}$ be the five Platonic solids: the tetrahedron, the cube, the octahedron, the dodecahedron and the icosahedron, respectively.
We also let $\xi_1(C_8)$, $\xi_2(C_8) $, $\xi_1(D_{20})$ and $\xi_2(D_{20})$ be polytopes shown in Figure~\ref{fig:dodeca}, and $B_n$ the bipyramid with $n$ vertices.
\begin{theorem}[\cite{ch-ki:arxiv}]
  Let $P$ be a three-dimensional simplicial polytope. If $P$ is reducible and
  combinatorially rigid, then $P$ is either $T_4 \# T_4 \# T_4$ or $P_1\#P_2$,
  where
\begin{align*}
  P_1 & \in \{T_4,O_6,I_{12}\},\\
  P_2 &\in \{T_4,O_6,I_{12},\xi_1(C_8),\xi_2(C_8)
  ,\xi_1(D_{20}),\xi_2(D_{20})\} \cup \{B_n:n\geq7\}.
\end{align*}
In addition, the following are combinatorially rigid in $\K_{poly}$:
$$T_4\# T_4, T_4\# O_6, T_4\# I_{12}, T_4\# B_n, O_6\# O_6, O_6\# B_n,$$
where $n\geq 7$. See Table~\ref{tab:red_rigid}.
\end{theorem}

\begin{table}
  \centering
  \begin{tabular}{c|c|c|c|c|c|c|c|c}
    $\#$ & $T_4$ & $O_6$ & $I_{12}$ & $B_n$, $n\geq7$ & $\xi_1(C_8)$& $\xi_2(C_8)$&
    $\xi_1(D_{20})$ & $\xi_2(D_{20})$\\ \hline
    $T_4$ & rigid & rigid & rigid & rigid & ?& ?& ?& ?\\ \hline
    $O_6$ & - & rigid & ? & rigid & ?& ?& ?& ?\\ \hline
    $I_{12}$ & - & - & ? & ? & ?& ?& ?& ?\\
\end{tabular}
  \caption{Reducible rigid polytopes of dimension $3$}\label{tab:red_rigid}
\end{table}

\end{document}